\documentclass[11pt,psfig,epsfig]{article}
\usepackage{amssymb,latexsym}
\usepackage{amsmath}
\usepackage{graphicx}
\newtheorem{theo}{Theorem}[section] %
\newtheorem{pro}[theo]{Proposition}

\begin{document}
\bibliographystyle{ieeetr}
\title{Segmented Tau Approximation for a Non-Autonomous Functional Differential Equation of Mixed Type}
\author{Carmen Da Silva\thanks{Escuela de Matem\'atica, Facultad de Ciencias,
Universidad Central de Venezuela, Caracas 1020, Venezuela \hspace{1cm}
 (carmen.dasilva@ciens.ucv.ve).}
 \thanks{This author was partially supported by the Consejo de Desarrollo Cient\'{\i}fico y Human\'{\i}stico (CDCH) at UCV.}
     \and
Ren\'e Escalante \thanks{Departamento de C\'omputo Cient\'{\i}fico y Estad\'{\i}stica,
 Divisi\'on de Ciencias F\'{\i}sicas y Matem\'aticas, Universidad Sim\'on
Bol\'{\i}var, Ap. 89000, Caracas, 1080-A, Venezuela
(rescalante@usb.ve).}
 \thanks{Corresponding author.}
\thanks{This author was partially supported by the Decanato de Investigaci\'on y Desarrollo (DID) at USB.}}
 \date{August 3, 2016}
\maketitle
\begin{abstract}
The segmented formulation of the Tau method is used to numerically solve the non-autonomous
forward-backward functional differential equation
\begin{equation*}
  \dot{x}(t)=a(t)x(t)+b(t)x(t-1)+c(t)x(t+1),
\end{equation*}
where $x$ is the unknown function, $a$, $b$, and $c$ are known functions. The step by step Tau method is applied to approximate the solution of this equation by a piecewise polynomial function.  A boundary value problem is posed, numerically solved, and analyzed. 
Also, a novel way to generate a set of non-autonomous problems with known analytical solution is provided.
From it, several non-autonomous problems were constructed and resolved with the proposed method.
We conclude that the good numerical results obtained in our numerical experimentation and the relative simplicity of
the Tau method demonstrate that it is a promising strategy
 for numerically solving mixed-type problems, as presented here. \\ [2mm] 
\end{abstract}
\newpage

\section{Introduction and Preliminaries}
In this paper the segmented Lanczos-Tau method is used to find numerical solutions for a family of
functional differential problems with both delayed and advanced arguments (i.e., forward-backward problems),
which are also referred to as mixed-type functional differential equations (MFDEs).
Here the following non-autonomous MFDE is considered
\begin{equation}\label{MTFDE}
    \dot{x}(t)=a(t)x(t)+b(t)x(t-1)+c(t)x(t+1),
\end{equation}
where $x(t)$ is a real single-valued function, and the coefficients $a(t)$, $b(t)$, and $c(t)$ are
known 
functions at $\mathbb{R}$. When the coefficients are constant, we say that the equation (\ref{MTFDE}) is autonomous;
 if at least one of them, is not a constant function, then (\ref{MTFDE}) is a non-autonomous equation.

The study of such equations is relatively recent and is motivated by the interest of some researchers
to attack some practical problems that arise naturally in various contexts, such as: the
modeling of the propagation of nerve impulse in a myelinated axon \cite{ChiBH:86}, problems in optimal
control \cite{Rusticha:89}, in economic dynamic problems \cite{Rustichb:89}, and
 in the study of traveling waves in discrete spatial media such as lattices \cite{MalVer:04}, \cite{ChCJMPN:00}, \cite{WuZou:97}.
 Recently, the MFDEs have been used in the analysis and modelling of economic growth (see \cite{Bou:04}, \cite{dAlAug:09}).
 It is clear that, with increasing frequency, significant and new applications that lead to mathematical models
 involving such equations appear, which undoubtedly has increased the need and interest of the scientific community
 to study the MFDEs, from a theoretical and numerical standpoint.

 While it is true that many of these models lead to non-linear MFDEs or systems of equations of this type, it is natural
 to develop theoretical and computational tools considering a much simpler family of equations, such as (\ref{MTFDE}),
 and then extending these insights to non-linear and vectorial cases. In this regard, Rustichini in \cite{Rusticha:89}
 addressed the spectrum of the linear (unbounded) operator, and constructed continuous semigroups on the stable, center, and
 unstable subspaces. His work led him to study two aspects of the theory of non-linear MFDEs: Hopf bifurcation and the center
 manifold theorems \cite{Rustichb:89}. 
 Likewise, Mallet-Paret and Verduyn \cite{MalVer:04} showed that for autonomous equations, the set of all forward solutions
 defines a semigroup which can be realized by a retarded functional differential equation, and similarly for the set of
 backward solutions as an advanced functional differential equation, where holomorphic factorizations played an important
 role in their results. Also, in \cite{HaSaSc:04}, an extension of this type of factorization for non-autonomous case is considered.

 We consider here the equation (\ref{MTFDE}) defined in the interval $[-1,\mathcal{K}]$ with $\mathcal{K}\geq 2$ a positive integer,
 subject to the following boundary conditions:
 \begin{eqnarray}
    x(t) & = & \psi_{1}(t), \hspace{.3cm} t\in[-1,0],   \label{CB1}  \\
    x(t) & = & \psi_{2}(t), \hspace{.3cm} t\in(\mathcal{K}-1,\mathcal{K}], \label{CB2}
 \end{eqnarray}
 where $\psi_{1}$ and $\psi_{2}$ are given continuous
functions in $[-1,0]$ and $[\mathcal{K}-1,\mathcal{K}]$,
 respectively, and satisfying equation (\ref{MTFDE}).

 Researchers Lima, Teodoro, Ford, and Lumb \cite{TeoLFL:09,FoLuLiTe:10,LiTeFoLu:10,LTFL:10} developed numerical algorithms to estimate the solution of (\ref{MTFDE})-(\ref{CB2}) in $(0,\mathcal{K}-1]$. They used linear $\theta$-methods, centered finite differences, collocation, least squares, and finite element methods.

In \cite{DaSEsc:11}, the step by step version of the Tau method was proposed by the authors to solve numerically
the boundary value problem (\ref{MTFDE})-(\ref{CB2}) in the autonomous case.
There we observe that our numerical results were consistent with those reported by other authors using other numerical approaches
($\theta$-method, least squares, and collocation methods, i.e. the above-mentioned methods), and where we obtained very satisfactory results.
We think that, in a similar way as in \cite{DaSEsc:11} the segmented Tau method was applied to the autonomous case,
it can also be extended to more complicated non-autonomous problem (\ref{MTFDE})-(\ref{CB2}).
By applying this method we seek to approximate the solution of equation (\ref{MTFDE}) by a piecewise polynomial function \cite{Ortiz:75}.
 In \cite{DaSEsc:14} a first approach, using the Tau method, was briefly presented to estimation
 of the solution of problem (\ref{MTFDE})-(\ref{CB2}).

\subsubsection*{The Tau method.}
 The {\it Tau method}, first introduced by Lanczos (\cite{Lanczos:38}, \cite{Lanczos:52} and \cite{Lanczos:56}),
 is based on the idea of getting efficient approximations of functions implicitly defined by a differential equation.
 This method, by construction, allows us to directly obtain polynomial approximations of high accuracy and reliability \cite{Ortiz:69}.
 It had been also reported in \cite{OrtPurRod:72} that the original
Tau method is, in many cases, comparable to the accuracy of best uniform
approximations or near optimal polynomial approximations of the same degree.

 The basic philosophy of the Tau method was extended to the numerical solution of linear and non-linear initial value, boundary value, and mixed problems for ordinary differential equations (see \cite{OnuOrt:84}, \cite{OrtPha:85} and \cite{OrtSam:81}), consequently applied to the eigenvalue problems (\cite{ChaOrt:68}, \cite{LiuOrt:87}, \cite{OrtSam:83}), to ``stiff" problems \cite{OnuOrt:84},  and to partial differential equations \cite{NamOrt:85}, among others. The Tau method has also been used as an analytic tool in the discussion of equivalence results across numerical methods (\cite{ElOrt:98}, \cite{OrtPha:97}).
It is an important feature of the Tau method that no trial solutions, approximate quadratures or large matrix inversions are required \cite{OnuOrt:84}.

A convergence analysis and error bounds for the Tau method was considered by Lanczos \cite{Lanczos:38}, \cite{Lanczos:73}, Luke \cite{Luke:69}, Ortiz and Pham Ngoc Dinh \cite{OrtPha:85}, \cite{OrtPha:84}, and El-Daou and Ortiz \cite{ElOrt:93}. The {\it recursive form of the Tau method}, formalized by Ortiz in \cite{Ortiz:69}, was extended to the case of systems of ordinary differential equations in \cite{FreOrt:82}, and also an error analysis was given there.
 In \cite{RooPfe:89} Roos and Pfeifer showed that
 the Tau method is, in the most interesting cases, a method of Galerkin-Petrov
 type, thus the convergence of the
 method follows from Vainikko's convergence theorem \cite{Vainikko:76}.

In the formulation of a {\it step by step} Tau version it is
allowed to construct {\it piecewise polynomial approximations} of a given function which can be used to start a refining process (see \cite{Ortiz:75} for details). The {\it Ortiz Step by Step Tau method} (or SST method to abbreviate) was later applied efficiently to the solution of linear and nonlinear boundary value problems \cite{OnuOrt:84}.
Computational strategies for a parallel implementation of the SST method were proposed in \cite{Escalan:03}.

In the Tau method, a {\it perturbation term} is introduced into the differential equation and from the perturbed equation an exact polynomial solution is obtained (Tau solution). This solution is an approximation to the solution of the original differential equation. In the segmented version of the Tau method (the SST method), the interval under consideration is divided into subintervals and the Tau method is applied separately in each subinterval. The Tau solution obtained in one interval is used as an input in the next interval. In \cite{CorEsca:07}, \cite{CorEscb:07}, and \cite{KhaOrt:96} the differential equation in each one of the subintervals is shifted to a corresponding equation in the interval $[0,1]$ (we apply this strategy here). This way, a sequence of differential equations defined in $[0,1]$ is established with the Tau solution for each equation providing information to its successor.

 In papers \cite{CorEsca:07}, \cite{CorEscb:07}, \cite{Khajah:05}, and \cite{KhaOrt:96} the step by step Tau method was applied to find polynomial approximations to the solution of linear and nonlinear delay differential equations. Also recently in \cite{DaSEsc:11} the segmented Tau method was applied to find approximations to the solution of an autonomous mixed-type functional differential boundary value problem.
 These papers seem to show that the segmented Tau method is a natural and promising strategy in the numerical solution of functional differential equations.

 For a brief exposition of the recursive formulation of the Tau method and the denominated {\it canonical polynomials},
 the reader is refereed to \cite{DaSEsc:11}, \cite{CorEsca:07}  and the references cited therein.

    This paper is organized as follows.
    In Section~2, we find the piecewise polynomial approximation of the non-autonomous boundary value problem (\ref{MTFDE})-(\ref{CB2})
    using the recursive formulation of the segmented Lanczos-–Tau method.
    In Section~3 we provide a novel way to generate a set of non-autonomous problems with known analytical solution,
    from which we extract numerical examples 
    to carry out comparisons.
    Finally, in Section~4, we present some concluding remarks.

\section{Solving non-autonomous problem using the segmented Tau method}
Our goal is to find a piecewise polynomial function, which 
is reasonably approximate to the exact solution of non-autonomous equation (\ref{MTFDE}) in the interval
$(0,\mathcal{K}-1]$, using information provided by the boundary conditions (\ref{CB1})-(\ref{CB2}).
 We assume that $a(t)$, $b(t)$, $c(t)$, $\psi_{1}(t)$, and $\psi_{2}(t)$ are polynomial functions
(otherwise, we will consider an accurate polynomial approximation associated with each of these functions).
Starting with the interval $(0,\mathcal{K}-1]$, it is partitioned into $\mathcal{K}-1$ unit subintervals,
denoted by $(k,k+1]$ with $k=0,1,\ldots ,\mathcal{K}-2$. The $k$-th subinterval, $(k,k+1]$, is associated
with step $k$ in the process of applying the segmented Tau method.

Let $s=t-k$; if $t\in(k,k+1]$ then $s\in[0,1]$. With this scaling strategy, we define:
\begin{equation*}
    x_{k}(s):=x(s+k)=x(t); \hspace{0.3cm} k=0,1,\ldots,\mathcal{K}-2 \hspace{0.3cm} \text{and} \hspace{0.3cm} s\in[0,1].
\end{equation*}
Thus, the problem (\ref{MTFDE})-(\ref{CB2}) is expressed as follows:
\begin{equation}\label{ProbDiscret}
    \left\{
      \begin{array}{ll}
        \dot{x}_{k}(s)=a_{k}(s)x_{k}(s)+b_{k}(s)x_{k-1}(s)+c_{k}(s)x_{k+1}(s), \\
        \vspace{-.4cm} \\
        x_{-1}(s)=\psi_{1}(s-1), \\
        x_{\mathcal{K}-1}(s)=\psi_{2}(s+\mathcal{K}-1),
      \end{array}
    \right. \vspace{.1cm}
\end{equation}
for $k=0,\ldots,\mathcal{K}-2$ and $s\in[0,1]$, where $a_{k}(s)$, $b_{k}(s)$, and $c_{k}(s)$ are
the corresponding representations of polynomials $a(t)$, $b(t)$, and $c(t)$, $t\in(k,k+1]$,
in the interval $[0,1]$. If we denote the degrees of the polynomials $a_{k}(s)$, $b_{k}(s)$, and $c_{k}(s)$ by
$d_{a}$, $d_{b}$, and $d_{c}$, respectively, then we can represent them by
\begin{equation}\label{Polinomiosabc}
    a_{k}(s)=\sum\limits_{i=0}^{d_{a}}\alpha_{i}^{(k)}s^{i}, \hspace{0.7cm} b_{k}(s)=\sum\limits_{i=0}^{d_{b}}\beta_{i}^{(k)}s^{i}, \hspace{0.7cm} c_{k}(s)=\sum\limits_{i=0}^{d_{c}}\gamma_{i}^{(k)}s^{i},
\end{equation}
where coefficient superscript $(k)$ indicates that it is currently taking into account the $k$-th step.

Let $D^{(k)}$ be the linear differential operator associated with (\ref{ProbDiscret}). I.e.,
\begin{equation}\label{OperadorDif}
    D^{(k)}[x_{k}(s)]:=\frac{d}{ds}[x_{k}(s)]-a_{k}(s)[x_{k}(s)], \hspace{0.6cm} k=0,\ldots,\mathcal{K}-2.
\end{equation}
As $k$ runs from $0$ to $\mathcal{K}-2$ in (\ref{OperadorDif}), we obtain $\mathcal{K}-1$ linear differential
operators which differ from each other by the polynomial $a_{k}(s)$, $s\in[0,1]$.
Next, using the Tau method, we will focus on finding a (piecewise) polynomial approximation of the solution of 
 (\ref{ProbDiscret}), so we must study how will be the set of canonical polynomials associated with the differential
operator (\ref{OperadorDif}) for fixed $k$.

\subsection{Canonical polynomials}
Lanczos proposed in 1956 \cite{Lanczos:56} the concept of the canonical polynomials $Q_n$,
$n \in \mathbb{N}_0\equiv \mathbb{N} \cup \{0\}$, associated with a linear differential operator.
Afterwards, Ortiz introduced the more workable definition of $Q_n$ \cite{Ortiz:69}:
\begin{equation}\label{DefPolCanOrtiz}
    D\left[Q_{n}(t)\right]=t^{n}+R_{n}(t), \hspace{.5cm} n\in\mathbb{N}_0-\mathcal{S},
\end{equation}
where $R_n(t)$ is a polynomial generated by $\{t^i\}$, $i\in S$ (set of indices for which canonical polynomials remain
undefined), and is called the {\it residual polynomial} of $Q_n(t)$. Another related concept is that of {\it generating polynomials} \cite{Ortiz:69}, which are obtained
from applying the associated differential operator to a power of $t$ (from them we can find a recursive relation
for the canonical polynomials).

For the case of the operator (\ref{OperadorDif}), the generating polynomial is given by
\begin{equation*}
    P_{m}^{(k)}(s):=D^{(k)}\left[s^{m}\right]=ms^{m-1}-\sum\limits_{i=0}^{d_{a}}\alpha_{i}^{(k)}s^{m+i},
\end{equation*}
whose degree is $\sigma_{m}=m+d_{a}$. 
Using the theory developed by Ortiz in \cite{Ortiz:69} we can prove the following result.

\begin{pro}\label{Pro-PolCanInd}
Let $k\in\{0, 1,\ldots,\mathcal{K}-2\}$ with $\mathcal{K}\geq 2$. The linear differential operator (\ref{OperadorDif})
 has $d_{a}$ undefined canonical polynomials.
\end{pro}

At this point, however, it would be helpful to consider the definition originally proposed by Lanczos for canonical polynomials \cite{Lanczos:56}:
\begin{equation}\label{DefPolCanLanczos}
    D^{(k)}\left[Q_{m}^{(k)}(s)\right]=s^{m} \hspace{.5cm} \text{with} \hspace{.3cm} s\in[0,1].
\end{equation}
From (\ref{DefPolCanLanczos}) and due to the linearity of the operator $D^{(k)}$ in the definition of the generating polynomial,
we obtain that
 $s^{m}=m\,Q_{m-1}^{(k)}(s)-\sum_{i=0}^{d_{a}}\alpha_{i}^{(k)}Q_{m+i}^{(k)}(s)$,
from which it follows that
\begin{equation}\label{FormRecur}
    Q_{m+d_{a}}^{(k)}(s)=-\frac{1}{\alpha_{d_{a}}^{(k)}}\left[s^{m}-m\,Q_{m-1}^{(k)}(s)+\sum\limits_{i=0}^{d_{a}-1}\alpha_{i}^{(k)}Q_{m+i}^{(k)}(s)\right], \hspace{.4cm} m\geq 0.\\
\end{equation}
When $0\leq m\leq d_{a}$, formula (\ref{FormRecur}) contains terms wherein the undefined canonical polynomials
$Q_{0}^{(k)}(s)$, $Q_{1}^{(k)}(s)$, $ \ldots, Q_{d_{a}-1}^{(k)}(s)$ appear. Proposition \ref{Pro-PolCanYRes}
shows a recursive formula that involves only polynomials whose indices are in $\mathbb{N}_0-S$ and
also provides the residual polynomials associated with these canonical polynomials.

\begin{pro}\label{Pro-PolCanYRes}
Let $d_{a}\neq 0$. If for each $m\geq 0$, $Q_{d_{a}+m}^{(k)}(s)$ is given by (\ref{FormRecur}), then
\begin{equation}\label{PolCan}
\begin{array}{c}
  Q_{d_{a}+m}^{(k)}(s)=\left\{
                     \begin{array}{lll}
                       -\frac{1}{\alpha_{d_{a}}^{(k)}}\left(s^{m}+\sum\limits_{i=0}^{m-1}\alpha_{d_{a}-m+i}^{(k)}Q_{d_{a}+i}^{(k)}(s)\right) \hspace{-.1cm} &\hspace{-.1cm}\text{if}& 0\leq m\leq d_{a}, \\
                     -\frac{1}{\alpha_{d_{a}}^{(k)}}\left(s^{m}-m\,Q_{m-1}^{(k)}(s)+\sum\limits_{i=1}^{d_{a}}\alpha_{d_{a}-i}^{(k)}Q_{d_{a}+m-i}^{(k)}(s)\right) \hspace{-.1cm} &\hspace{-.1cm}\text{if}& m\geq d_{a}+1,
                     \end{array}
                   \right.
\end{array}
\end{equation}
are the canonical polynomials associated with $D^{(k)}$ ($k=0,\ldots,\mathcal{K}-2$), and
\begin{equation}\label{PolRes}
    R_{d_{a}+m}^{(k)}(s)=\left\{
                     \begin{array}{lll}
                       \frac{1}{\alpha_{d_{a}}^{(k)}}\left(-m\,s^{m-1}+\sum\limits_{i=0}^{d_{a}-m-1}\alpha_{i}^{(k)}s^{m+i}\right) &\text{if}& 0\leq m\leq d_{a}, \\
                       0 &\text{if}& m\geq d_{a}+1,
                     \end{array}
                   \right.
\end{equation}
are the corresponding residual polynomials.
\end{pro}

\noindent{\it Note}:
Here we make the convention that if $m=0$, the sum of the first part of $Q_{d_{a}+m}^{(k)}(s)$ is equal to zero.
 Similarly, when $m=d_{a}$ in $R_{d_{a}+m}^{(k)}(s)$.

\noindent{\it \textbf{Proof.}} When $m\in\{0,\ldots,d_{a}\}$, undefined canonical polynomials correspond to the
 indices: $m$, $m+1$, $\ldots$, $m+(d_{a}-m-1)$. From (\ref{FormRecur}) we obtain the canonical polynomial of index $d_{a}+m$,
 neglecting terms involving an undefined canonical polynomial, since the coefficients of these terms become part
 of the residual polynomials. Thus, we obtain the first part of (\ref{PolCan}).
 When $0\leq m\leq d_{a}$, the residual polynomial associated with $Q_{d_{a}+m}^{(k)}(s)$
 can be derived by applying the operator $D^{(k)}$ to both sides of (\ref{PolCan});
 then, from  (\ref{OperadorDif}) and definition (\ref{DefPolCanLanczos}) it follows that
\begin{equation*}
    D^{(k)}\left[Q_{d_{a}+m}^{(k)}(s)\right]=-\frac{1}{\alpha_{d_{a}}^{(k)}}\left(m\,s^{m-1}-\sum\limits_{i=0}^{d_{a}-m-1}\alpha_{i}^{(k)}s^{m+i}\right)+s^{d_{a}+m}.
\end{equation*}
We now replace the left side of the above equation by the appropriate expression of definition (\ref{DefPolCanOrtiz})
and solve for $R_{d_{a}+m}^{(k)}(s)$ to obtain the first part of (\ref{PolRes}).\\

In the case that $m\in\{d_{a}+1,d_{a}+2,\ldots\}$, we deduce from formula (\ref{FormRecur}) that there is no undefined canonical polynomials.
From this it follows that the second part of (\ref{PolCan}) is a rewriting of (\ref{FormRecur}).
Also, if we perform a procedure similar to the one described in the previous paragraph to attain the residual polynomials,
 we obtain in this case that $R_{d_{a}+m}^{(k)}(s)=0$.
 $\Box $

\begin{pro}
    For $m\geq 0$, the canonical polynomial $Q_{d_{a}+m}^{(k)}(s)$, defined recursively by (\ref{PolCan}), is a polynomial of degree $m$.
\end{pro}

\noindent{\it \textbf{Proof.}}
    It follows by induction on $m$.  $\Box $

We now introduce the following notation. For $k=0.\ldots,\mathcal{K}-2$,
\begin{equation}\label{NotacPolCan}
    Q_{d_{a}+m}^{(k)}(s)=\sum\limits_{j=0}^{m}q_{j}^{(k)(m)}s^{j}, 
\end{equation}
where $q_{j}^{(k)(m)}$ represents the $j$th coefficient of the polynomial $Q_{d_{a}+m}^{(k)}(s)$, which can be computed from (\ref{PolCan}).
Note that as $m\geq 0$, for fixed $k$, we have an infinite number of canonical polynomials.
 Therefore, the superscript $(m)$ on $q_{j}^{(k)(m)}$ refers to the $m$th polynomial associated with the operator $D^{(k)}$.

\subsection{Solution of perturbed non-autonomous problem} \label{sectsolpert}
The Tau method, as originally proposed by Lanczos \cite{Lanczos:52}, is based on disturbing the differential equation
$D[x(t)]=0$, $t\in J$ (a closed interval in $\mathbb{R}$),  by putting an error term on the right side of the equation.
This term is introduced intentionally in order to make the equation solvable by a finite power series.
The perturbed differential equation is given by $D[x_{n}(t)]=H(t)$, $t\in J$, where $H(t)=\left(\sum_{i=0}^{r}\tau_{i}t^{i}\right)T_{n-r}(t)$, which is called the {\it perturbation term},
  $T_{n-r}$ is the Chebyshev polynomial of degree $n-r$ on $J$, and the $\tau_i$ parameters are chosen such that $x_{n}$
 satisfies exactly supplementary conditions, i.e. the initial, boundary or mixed conditions. So, $x_{n}$
 is a polynomial approximation to $x$ in $J$.

In our case, for each $k=0,\ldots,\mathcal{K}-2$, we perturb the differential equation (\ref{ProbDiscret}).
For this, we need $\mathcal{K}-1$ perturbation terms such that the exact solution of the $k$th perturbed differential equation
is a polynomial of degree $n$.
 The coefficients of the $\mathcal{K}-1$ polynomial solutions and $\tau$-parameters of perturbation terms
 will be determined by solving a square system of linear equations (see Theorem~\ref{Teo-SolEDperturbNoAut}).
 Therefore, we need to consider $d\left(\mathcal{K}-1\right)+\mathcal{K}$  $\tau$-parameters in all the steps of the process,
 where $d=\max\left\{d_{a},d_{b},d_{c}\right\}$ (i.e., $d+2$ $\tau$-parameters for step 0, and $d+1$ $\tau$-parameters
 for each of the remaining $\mathcal{K}-2$ steps).
 Accordingly, for $s\in[0,1]$, the perturbation term becomes
\vspace{.2cm}\begin{equation}\label{TermErrorProbNoAut}
    H_{d+n}^{(k)}(s)=\left\{
                      \begin{array}{ll}
                        \left(\sum\limits_{i=0}^{d+1}\tau_{i}^{(0)}s^{i}\right)T_{n-1}^{*}(s), & \hbox{$k=0$;} \\
                        \left(\sum\limits_{i=0}^{d}\tau_{i}^{(k)}s^{i}\right)T_{n}^{*}(s), & \hbox{$1\leq k\leq \mathcal{K}-2$,}
                      \end{array}
                    \right.\vspace{.2cm}
\end{equation}
where $T_{n}^{*}$ represents the shifted Chebyshev polynomial of degree $n$ in $[0,1]$ (it is similar for $T_{n-1}^{*}$).

Then, the perturbed non-autonomous problem, for $k=0,\ldots,\mathcal{K}-2$ with $\mathcal{K}\geq 2$, becomes
\begin{equation}\label{ProblemaPerturbado}
    \left\{
      \begin{array}{ll}
        D^{(k)}[X_{k}(s)]=b_{k}(s)X_{k-1}(s)+c_{k}(s)X_{k+1}(s)+H_{d+n}^{(k)}(s); \hspace{.4cm} s\in[0,1], \\
        \vspace{-.4cm} \\
        X_{-1}(s)=\psi_{1}(s-1), \\
        X_{\mathcal{K}-1}(s)=\psi_{2}(s+\mathcal{K}-1),
      \end{array}
    \right. \vspace{.1cm}
\end{equation}
 where $D^{(k)}$ is defined as in (\ref{OperadorDif}).\\

The following theorem provides information about solving the problem (\ref{ProblemaPerturbado}) when $a_{k}(s)$, $b_{k}(s)$, $c_{k}(s)$, $\psi_{1}$, and $\psi_{2}$ are polynomials.

\begin{theo}\label{Teo-SolEDperturbNoAut}
Consider the problem (\ref{ProblemaPerturbado}) with $a_{k}$, $b_{k}$, $c_{k}$, $\psi_{1}$, and $\psi_{2}$ polynomial functions
 of degree $\leq n$. For every $k=0,\ldots,\mathcal{K}-2$, with $\mathcal{K}\geq 2$, the exact solution of the problem (\ref{ProblemaPerturbado}) is determined by a polynomial of degree $n$ of the form
\begin{equation}\label{XsolEDperturb}
    X_{k}(s)=\sum\limits_{i=0}^{n}a_{i}^{(k)}s^{i}, \hspace{.6cm} s\in[0,1]. \\
\end{equation}
Moreover, if $d=d_{a}$, $d_{a}\neq 0$, and
\begin{equation}\label{XCondContin}
    X_{k}(0)=X_{k-1}(1) \hspace{0.4cm} \text{with} \hspace{0.4cm} k=0,\ldots,\mathcal{K}-1,
\end{equation}
then for every step $k$, there exists a square system of linear equations whose solution provides the polynomial equation coefficients   $\left\{a_{i}^{(k)}\right\}_{i=0}^{n}$, and perturbation term $\tau$-parameters $\left\{\tau_{i}^{(0)}\right\}_{i=0}^{d+1}$ and $\left\{\tau_{i}^{(k)}\right\}_{i=0}^{d}$.
\end{theo}
\noindent{\it Note}:
The $\mathcal{K}$ conditions arising from (\ref{XCondContin}) refer to the continuity of the polynomial
 solution of perturbed non-autonomous differential equation, through successive steps;
 also providing a continuous connection to the boundary conditions.

\noindent{\it \textbf{Proof.}}
Our goal will be to define a system of the form $A\vec{x}=\vec{b}$, where the unknown vector $\vec{x}$
 can have its elements stored in blocks, as follows:
\begin{equation*}
    \vec{x}=\left(
       \begin{array}{cccccccc}
         \vec{a}^{(0)},
         \vec{a}^{(1)},
         \ldots,
         \vec{a}^{(\mathcal{K}-2)},
         \vec{\tau}^{(0)},
         \vec{\tau}^{(1)},
         \ldots,
         \vec{\tau}^{(\mathcal{K}-2)}
       \end{array}
     \right)^{T},
\end{equation*}
 where, for $k=0,1,\ldots,\mathcal{K}-2$, $\vec{a}^{(k)}=\left(a_{0}^{(k)}, a_{1}^{(k)}, \ldots, a_{n}^{(k)}\right)^{T}\in\mathbb{R}^{n+1}$,
  $\vec{\tau}^{(0)}=\left(\tau_{0}^{(0)}, \tau_{1}^{(0)}, \ldots, \tau_{d}^{(0)}, \tau_{d+1}^{(0)}\right)^{T}\in\mathbb{R}^{d+2}$, and, for $k=1,2,\ldots,\mathcal{K}-2$, $\vec{\tau}^{(k)}=\left(\tau_{0}^{(k)}, \tau_{1}^{(k)}, \ldots, \tau_{d}^{(k)}\right)^{T}\in\mathbb{R}^{d+1}$. $T$ denotes the vector transpose.

Next, we define $a_{i}^{(-1)}$ and $a_{i}^{(\mathcal{K}-1)}$, $i=0,\ldots,n$, as coefficients of polynomial boundary conditions
$X_{-1}$ and $X_{\mathcal{K}-1}$, respectively, according to the notation introduced in (\ref{XsolEDperturb}).
 If for any boundary condition, the degree, say $\nu$, is less than $n$, we make the convention that $a_{i}^{(-1)}=0$ or $a_{i}^{(K-1)}=0$
  for all $i=\nu+1,\ldots,n$. In addition, we will use the notation $C_{j}^{(n-1)}$, $j=0,1,\ldots,n-1$, and $C_{j}^{(n)}$, $j=0,1,\ldots,n$,
 in (\ref{TermErrorProbNoAut}) to refer to the coeficients of the shifted Chebyshev polynomials of degrees $n-1$ and $n$, respectively.

Taking into account that $d=d_{a}\in\left\{1,\ldots,n\right\}$, four possible cases can be considered.
\begin{flushleft}
\textsc{Case 1:} $d=d_{a}=d_{b}=d_{c}$.
\end{flushleft}
Since the perturbation term is a function defined piecewise, we also need to study two cases: $k=0$ and $k\in\{1,2,\ldots,\mathcal{K}-2\}$.

\noindent{\it Case} $k=0$. If we substitute (\ref{TermErrorProbNoAut}) for $k=0$, (\ref{XsolEDperturb}) with $k=-1$ and $k=1$, and
 polynomials $b_{0}(s)$ and $c_{0}(s)$ using (\ref{Polinomiosabc}), on the right side of the differential equation
 (\ref{ProblemaPerturbado}) (with $k=0$), and if we rewrite the obtained expression in a more suitable form we have
\begin{equation*}
     \begin{array}{l}
       D^{(0)}[X_{0}(s)] = \sum\limits_{i=0}^{d}\left(\sum\limits_{j=0}^{i}\beta_{j}^{(0)}a_{i-j}^{(-1)}\right)s^{i} + \sum\limits_{i=1}^{n-d}\left(\sum\limits_{j=0}^{d}\beta_{j}^{(0)}a_{d+i-j}^{(-1)}\right)s^{d+i} \\
        + \sum\limits_{i=1}^{d}\left(\sum\limits_{j=n-d}^{n-i}\beta_{n-j}^{(0)}a_{i+j}^{(-1)}\right)s^{n+i}
        + \sum\limits_{i=0}^{d}\left(\sum\limits_{j=0}^{i}\gamma_{j}^{(0)}a_{i-j}^{(1)}\right)s^{i} \\ + \sum\limits_{i=1}^{n-d}\left(\sum\limits_{j=0}^{d}\gamma_{j}^{(0)}a_{d+i-j}^{(1)}\right)s^{d+i}
        + \sum\limits_{i=1}^{d}\left(\sum\limits_{j=n-d}^{n-i}\gamma_{n-j}^{(0)}a_{i+j}^{(1)}\right)s^{n+i} \\
        + \sum\limits_{i=0}^{d-1}\left(\sum\limits_{j=0}^{i}\tau_{j}^{(0)}C_{i-j}^{(n-1)}\right)s^{i} + \sum\limits_{i=1}^{n-d}\left(\sum\limits_{j=0}^{d-1}\tau_{j}^{(0)}C_{d-1+i-j}^{(n-1)}\right)s^{d-1+i} \\
        + \sum\limits_{i=1}^{d-1}\left(\sum\limits_{j=n-d}^{n-1-i}\tau_{n-1-j}^{(0)}C_{i+j}^{(n-1)}\right)s^{n-1+i}
        + \sum\limits_{i=0}^{n-1}\tau_{d}^{(0)}C_{i}^{(n-1)}s^{d+i} + \sum\limits_{i=0}^{n-1}\tau_{d+1}^{(0)}C_{i}^{(n-1)}s^{d+1+i}.
     \end{array} \vspace{.2cm}
\end{equation*}
 If $d=n$, we make the convention that the sums with index of summation from $i=1$ to $n-d$ are zero.

By replacing $s^{i}$, $s^{d+i}$, $s^{n+i}$, $s^{d-1+i}$, $s^{n-1+i}$, and $s^{d+1+i}$ by its equivalent expression given in  (\ref{DefPolCanLanczos}), and then applying the linearity of the differential operator $D^{(0)}$,
 after a suitable grouping of terms, the following expression for $X_{0}$ is obtained,
\begin{equation}\label{primerX0}
     \begin{array}{l}
       X_{0}(s) = \sum\limits_{i=0}^{d}\left[\sum\limits_{j=0}^{i}\left(\beta_{j}^{(0)}a_{i-j}^{(-1)}+\gamma_{j}^{(0)}a_{i-j}^{(1)}\right)\right]Q_{i}^{(0)}(s)\\
                + \sum\limits_{i=1}^{n-d}\left[\sum\limits_{j=0}^{d}\left(\beta_{j}^{(0)}a_{d+i-j}^{(-1)}+\gamma_{j}^{(0)}a_{d+i-j}^{(1)}\right)\right]Q_{d+i}^{(0)}(s)\\
                + \sum\limits_{i=1}^{d}\left[\sum\limits_{j=n-d}^{n-i}\left(\beta_{n-j}^{(0)}a_{i+j}^{(-1)}+\gamma_{n-j}^{(0)}a_{i+j}^{(1)}\right)\right]Q_{n+i}^{(0)}(s)\\
                + \sum\limits_{i=0}^{d-1}\left(\sum\limits_{j=0}^{i}\tau_{j}^{(0)}C_{i-j}^{(n-1)}\right)Q_{i}^{(0)}(s)
       + \sum\limits_{i=1}^{n-d}\left(\sum\limits_{j=0}^{d-1}\tau_{j}^{(0)}C_{d-1+i-j}^{(n-1)}\right)Q_{d-1+i}^{(0)}(s)\\
                + \sum\limits_{i=1}^{d-1}\left(\sum\limits_{j=n-d}^{n-1-i}\tau_{n-1-j}^{(0)}C_{i+j}^{(n-1)}\right)Q_{n-1+i}^{(0)}(s)
       + \sum\limits_{i=0}^{n-1}\tau_{d}^{(0)}C_{i}^{(n-1)}Q_{d+i}^{(0)}(s)\\
                + \sum\limits_{i=0}^{n-1}\tau_{d+1}^{(0)}C_{i}^{(n-1)}Q_{d+1+i}^{(0)}(s).
     \end{array}
\end{equation}
As Proposition~\ref{Pro-PolCanInd} shows, the differential operator $D^{(0)}$ has $d$ undefined canonical polynomials
 corresponding to the indices $0,1,\ldots,d-1$, so it is assumed that the coefficients of $Q_{0}^{(0)}$, $Q_{1}^{(0)}$, $\ldots$, $Q_{d-1}^{(0)}$ in (\ref{primerX0}) are equal to zero. This leads to the following linear system of $d$ equations,
\begin{equation*}\label{X0coefcero}
    \sum\limits_{j=0}^{i}\left(\gamma_{j}^{(0)}a_{i-j}^{(1)}+\tau_{j}^{(0)}C_{i-j}^{(n-1)}\right)=-\sum\limits_{j=0}^{i}\beta_{j}^{(0)}a_{i-j}^{(-1)}, \hspace{.2cm} \text{for} \hspace{.2cm} i=0,\ldots,d-1. \vspace{.4cm} \\
\end{equation*}

Then, expression (\ref{primerX0}) reduces to
\begin{equation}\label{segundoX0}
    \begin{array}{l}
      X_{0}(s) = \sum\limits_{j=0}^{d}\left(\beta_{j}^{(0)}a_{d-j}^{(-1)}+\gamma_{j}^{(0)}a_{d-j}^{(1)}\right)Q_{d}^{(0)}(s)
      + \sum\limits_{i=1}^{n-d}\left[\sum\limits_{\ell=0}^{d}\left(\beta_{\ell}^{(0)}a_{d+i-\ell}^{(-1)}+\gamma_{\ell}^{(0)}a_{d+i-\ell}^{(1)}\right)\right]Q_{d+i}^{(0)}(s) \\
               + \sum\limits_{i=1}^{d}\left[\sum\limits_{\ell=n-d}^{n-i}\left(\beta_{n-\ell}^{(0)}a_{i+\ell}^{(-1)}+\gamma_{n-\ell}^{(0)}a_{i+\ell}^{(1)}\right)\right]Q_{n+i}^{(0)}(s)
      + \sum\limits_{i=1}^{n-d}\left(\sum\limits_{\ell=0}^{d-1}\tau_{\ell}^{(0)}C_{d-1+i-\ell}^{(n-1)}\right)Q_{d+i-1}^{(0)}(s) \\
               + \sum\limits_{i=1}^{d-1}\left(\sum\limits_{\ell=n-d}^{n-1-i}\tau_{n-1-\ell}^{(0)}C_{i+\ell}^{(n-1)}\right)Q_{n-1+i}^{(0)}(s)
      + \sum\limits_{i=0}^{n-1}\tau_{d}^{(0)}C_{i}^{(n-1)}Q_{d+i}^{(0)}(s) \\
               + \sum\limits_{i=0}^{n-1}\tau_{d+1}^{(0)}C_{i}^{(n-1)}Q_{d+1+i}^{(0)}(s).
    \end{array}
\end{equation}
Next, we replace (\ref{XsolEDperturb}), with $k=0$, on the left side of (\ref{segundoX0}) and apply the notation (\ref{NotacPolCan})
for the canonical polynomials.  All to conveniently rewrite (\ref{segundoX0}), equate coefficients on both
sides of the expression, and obtain $n+1$ additional linear equations.

\noindent{\it Case} $k\in\{1,\ldots,\mathcal{K}-2\}$. The procedure is analogous to that performed in the previous case.
Substituting into the differential equation (\ref{ProblemaPerturbado}) the expressions given in (\ref{Polinomiosabc})
 for $b_{k}$(s) and $c_{k}(s)$, the formula (\ref{XsolEDperturb}) with indices $k-1$ and $k+1$, and the second part
 of the perturbation term (\ref{TermErrorProbNoAut}), a polynomial expression to  $D^{(k)}\left[X_{k}(s)\right]$
 will be obtained. Finally, from (\ref{DefPolCanLanczos}) and the linearity of the differential operator $D^{(k)}$,
 it follows that
\begin{equation}\label{primerXk}
    \begin{array}{c}
      X_{k}(s)=\sum\limits_{i=0}^{d}\left[\sum\limits_{j=0}^{i}\left(\beta_{j}^{(k)}a_{i-j}^{(k-1)} + \gamma_{j}^{(k)}a_{i-j}^{(k+1)} + \tau_{j}^{(k)}C_{i-j}^{(n)}\right)\right]Q_{i}^{(k)}(s) \\
      + \sum\limits_{i=1}^{n-d}\left[\sum\limits_{j=0}^{d}\left(\beta_{j}^{(k)}a_{d+i-j}^{(k-1)} + \gamma_{j}^{(k)}a_{d+i-j}^{(k+1)} + \tau_{j}^{(k)}C_{d+i-j}^{(n)}\right)\right]Q_{d+i}^{(k)}(s) \\
      + \sum\limits_{i=1}^{d}\left[\sum\limits_{j=n-d}^{n-i}\left(\beta_{n-j}^{(k)}a_{i+j}^{(k-1)} + \gamma_{n-j}^{(k)}a_{i+j}^{(k+1)} + \tau_{n-j}^{(k)}C_{i+j}^{(n)}\right)\right]Q_{n+i}^{(k)}(s).
    \end{array}
\end{equation}

Note that Proposition \ref{Pro-PolCanInd} implies that equation (\ref{primerXk}) contains the undefined canonical polynomials
 $Q_{0}^{(k)}$, $Q_{1}^{(k)}$, $\ldots$, $Q_{d-1}^{(k)}$. Therefore, the coefficients of the terms involving these polynomials
 are equated to zero. That is,
\begin{equation}\label{Xkcoefcero}
    \sum\limits_{j=0}^{i}\left(\beta_{j}^{(k)}a_{i-j}^{(k-1)}+\gamma_{j}^{(k)}a_{i-j}^{(k+1)}+\tau_{j}^{(k)}C_{i-j}^{(n)}\right)=0 \hspace{.3cm} \text{for} \hspace{.3cm} i=0,\ldots,d-1.
\end{equation}
This leads to the simplification of (\ref{primerXk}) as follows 
\begin{equation}\label{segundoXk}
    \begin{array}{l}
      X_{k}(s)=\sum\limits_{j=0}^{d}\left(\beta_{j}^{(k)}a_{d-j}^{(k-1)} + \gamma_{j}^{(k)}a_{d-j}^{(k+1)} + \tau_{j}^{(k)}C_{d-j}^{(n)}\right)Q_{d}^{(k)}(s)\\
      + \sum\limits_{i=1}^{n-d}\left[\sum\limits_{\ell=0}^{d}\left(\beta_{\ell}^{(k)}a_{d+i-\ell}^{(k-1)}+\gamma_{\ell}^{(k)}a_{d+i-\ell}^{(k+1)}+\tau_{\ell}^{(k)}C_{d+i-\ell}^{(n)}\right)\right]Q_{d+i}^{(k)}(s)\\ + \sum\limits_{i=1}^{d}\left[\sum\limits_{\ell=n-d}^{n-i}\left(\beta_{n-\ell}^{(k)}a_{i+\ell}^{(k-1)}+\gamma_{n-\ell}^{(k)}a_{i+\ell}^{(k+1)}+\tau_{n-\ell}^{(k)}C_{i+\ell}^{(n)}\right)\right]Q_{n+i}^{(k)}(s).\\
    \end{array}
\end{equation}

Note that when $k=\mathcal{K}-2$, the expression $X_{k+1}(s)$ appearing in (\ref{ProblemaPerturbado}) becomes $X_{\mathcal{K}-1}(s)$, which represents one of the boundary conditions.
 By using the notation (\ref{NotacPolCan}) and (\ref{XsolEDperturb}) in (\ref{segundoXk}), equating the coefficients
 on both sides of this expression, we get a set of $n+1$ linear equations.
Furthermore, if rewrite (\ref{XCondContin}) using the notation (\ref{XsolEDperturb}), we obtain $\mathcal{K}$
additional linear equations. That is,
\begin{equation}\label{sistemacondcontin}
\left\{
  \begin{array}{l}
    a_{0}^{(k)}=\sum\limits_{i=0}^{n}a_{i}^{(-1)},\\
    \sum\limits_{i=0}^{n}a_{i}^{(k-1)}-a_{0}^{(k)}=0 \hspace{.4cm} \text{para} \hspace{.4cm} k=1,\ldots,\mathcal{K}-2,\\
    \sum\limits_{i=0}^{n}a_{i}^{(\mathcal{K}-2)}=a_{0}^{(\mathcal{K}-1)}.
  \end{array}
\right.
\end{equation}

So, we have a total of $(n+d+1)(\mathcal{K}-1)+\mathcal{K}$ linear equations, and the matrix (in blocks)
of coefficients of the system becomes,
\begin{equation}\label{matrizAprobNoAut}
    A\hspace{-.1cm}=\hspace{-.1cm}\left(
        \begin{array}{cccccccccc}
          I & U_{\gamma}^{(0)} &  &  &  & \widetilde{R}_{C} &  &  &  &  \\
          U_{\beta}^{(1)} & I & \ddots &  &  &  & R_{C}^{(1)} &  &  &  \\
            & \ddots & \ddots & \ddots &  &  &  & \ddots &  &  \\
           &  & \ddots & I & U_{\gamma}^{(\mathcal{K}-3)} &  &  &  & R_{C}^{(\mathcal{K}-3)} &  \\
           &  &  & U_{\beta}^{(\mathcal{K}-2)} & I &  &  &  &  & R_{C}^{(\mathcal{K}-2)} \\
          \textbf{0} & R_{\gamma}^{(0)} &  &  &  & \widehat{R}_{C} &  &  &  &  \\
          R_{\beta}^{(1)} & \textbf{0} & \ddots &  &  &  & R_{C} &  &  &  \\
           & \ddots & \ddots & \ddots &  &  &  & \ddots &  &  \\
           &  & \ddots & \textbf{0} & R_{\gamma}^{(\mathcal{K}-3)} &  &  &  & R_{C} &  \\
           &  &  & R_{\beta}^{(\mathcal{K}-2)} & \textbf{0} &  &  &  &  & R_{C} \\
          M_{1} & M_{2} & \ldots & M_{\mathcal{K}-2} & M_{\mathcal{K}-1} &  &  &  &  &  \\
        \end{array}
      \right),
\end{equation}
where $I$ represents the identity matrix of order $(n+1)\times (n+1)$. $U_{\gamma}^{(k)}$ and $U_{\beta}^{(k)}$
are upper triangular matrices with the same order as the matrix $I$. For $k=0,\ldots,\mathcal{K}-3$, non-zero elements
of the $j$th column of the matrix $U_{\gamma}^{(k)}$ are defined by
\begin{align*}
    &-\sum\limits_{\ell=0}^{d}\gamma_{d-\ell}^{(k)}\;q_{i-1}^{(k)(j-1-\ell)} \hspace{.7cm} \text{for} \hspace{.7cm} i=1,\ldots,j-d-1\\
    &-\sum\limits_{\ell=0}^{j-i}\gamma_{d-\ell}^{(k)}\;q_{i-1}^{(k)(j-1-\ell)} \hspace{.7cm} \text{for} \hspace{.7cm} i=j-d,\ldots,j.
\end{align*}
Similarly, for $k=1,\ldots,\mathcal{K}-2$, the non-zero elements of the $j$th column of $U_{\beta}^{(k)}$
are defined as those of $U_{\gamma}^{(k)}$ (exchanging $\gamma$ and $\beta$). The order of matrix
$\widetilde{R}_{C}$ is $(n+1)\times (d+2)$. When $d=n$, the first column of the matrix $\widetilde{R}_{C}$
has entries equal to zero, while the $ij$th non-zero element, $2\leq j\leq n+2$, is defined by 
\begin{equation*}
    -\sum\limits_{\ell=1}^{j-i}C_{n-\ell}^{(n-1)}q_{i-1}^{(0)(j-1-\ell)} \;\;\; \text{for} \;\;\; j>i,
\end{equation*}
where we make the convention that if $j=n+2$ and $i=1$ then $C_{-1}^{(n-1)}=0$.
When $d=n-1$, $\widetilde{R}_{C}$ is a square and upper triangular matrix, and its non-zero elements are as follows,
\begin{equation*}
    -\sum\limits_{\ell=0}^{j-i}C_{n-1-\ell}^{(n-1)}\;q_{i-1}^{(0)(j-1-\ell)} \;\;\; \text{for} \;\;\; j\geq i.
\end{equation*}
Now, if $d\leq n-2$, the non-zero elements of the $j$th column of the matrix $\widetilde{R}_{C}$
are as shown in the following expression,
\begin{equation*}
    -\sum\limits_{\ell=0}^{n-d+j-i-1}C_{n-1-\ell}^{(n-1)}\;q_{i-1}^{(0)(n-d+j-2-\ell)} \;\;\; \text{for} \;\;\; i=1,\ldots,n-d+j-1,
\end{equation*}
where $C_{-1}^{(n-1)}=0$ if $j=d+2$ and $i=1$. For $k=1,\ldots,\mathcal{K}-2$, the matrix $R_{C}^{(k)}$ is of order $(n+1)\times (d+1)$,
and its non-zero elements are of the form
\begin{equation*}
    -\sum\limits_{\ell=0}^{n-d+j-i}C_{n-\ell}^{(n)}\;q_{i-1}^{(k)(n-d+j-1-\ell)} \;\;\; \text{for} \;\;\; i=1,\ldots,n-d+j.
\end{equation*}
It is worth noting that if $d=n$, $R_{C}^{(k)}$ is a square and upper triangular matrix. $R_{\gamma}^{(k)}$ and $R_{\beta}^{(k)}$
are rectangular matrices of order $d\times (n+1)$ and their last $n-d+1$ columns are zero vectors. For $i=1,\ldots,j$,
the $ij$th non-zero element of $R_{\gamma}^{(k)}$, $k=0,\ldots,\mathcal{K}-3$, is $\gamma_{i-1}^{(k)}$,  and the
$ij$th non-zero element of the matrix $R_{\beta}^{(k)}$, $k=1,\ldots,\mathcal{K}-2$, is given by $\beta_{i-1}^{(k)}$.
Null matrices, $\textbf{0}$, with the same order as the above matrices are shown in $A$ to emphasize the separation
of the bands that contain $R_{\gamma}^{(k)}$ and $R_{\beta}^{(k)}$. The matrix $\widehat{R}_{C}$ is of order $d\times (d+2)$, and
 its $ij$th non-zero element is of the form $C_{i-1}^{(n-1)}$, for $i=1,\ldots,j$.
 Similarly, $R_{C}$ is a rectangular matrix of order $d\times (d+1)$ and its $ij$th non-zero element is $C_{i-1}^{(n)}$, for $i=1,\ldots,j$. In total, there are $K-2$ matrices $R_{C}$ in the $A$ matrix. Finally, the $M_{i}$ matrices, $i=1,\ldots,\mathcal{K}-1$, are of order $\mathcal{K}\times (n+1)$ and have two non-zero rows, namely:
  the row $i$, represented by $(1,0,\cdots,0)$, and the row $i+1$ of the form $(1,1,\cdots,1)$. So, the matrix $A$ given in
   (\ref{matrizAprobNoAut}) is of order $[(n+d+1)(\mathcal{K}-1)+\mathcal{K}]\times[(n+d+1)(\mathcal{K}-1)+\mathcal{K}]$.

The vector of independent terms $\vec{b}$ has a block structure as follows,
\begin{equation}\label{vectorbprobNoAut}
    \vec{b}=\left(
      \begin{array}{ccccccccccc}
        \vec{u}^{(-1)}, &
        \vec{\textbf{0}}, &
        \ldots, &
        \vec{\textbf{0}}, &
        \vec{u}^{(\mathcal{K}-1)}, &
        \vec{v}^{(-1)}, &
        \vec{\textbf{0}}, &
        \ldots, &
        \vec{\textbf{0}}, &
        \vec{v}^{(\mathcal{K}-1)}, &
        \vec{w}
      \end{array}
    \right)^{T},
    \vspace{.1cm}
\end{equation}
where $\vec{u}^{(-1)}$, $\vec{u}^{(\mathcal{K}-1)}$ $\in\mathbb{R}^{n+1}$,  $\vec{v}^{(-1)}$, $\vec{v}^{(\mathcal{K}-1)}$ $\in\mathbb{R}^{d}$, and $\vec{w}$ $\in\mathbb{R}^{\mathcal{K}}$. Null vectors located between $\vec{u}^{(-1)}$ and $\vec{u}^{(\mathcal{K}-1)}$ are of dimension $n+1$; however, those null vectors between $\vec{v}^{(-1)}$ and $\vec{v}^{(\mathcal{K}-1)}$ have dimension $d$. Thus, $\vec{b}\in\mathbb{R}^{(n+d+1)(\mathcal{K}-1)+\mathcal{K}}$.

The elements of $\vec{u}^{(-1)}$ are
\begin{equation*}
    \begin{split}
      u_{1}^{(-1)} &= \sum\limits_{\ell=0}^{d}\beta_{\ell}^{(0)}a_{d-\ell}^{(-1)}q_{0}^{(0)(0)} + \sum\limits_{j=1}^{n-d}\left(\sum\limits_{\ell=0}^{d}\beta_{\ell}^{(0)}a_{d+j-\ell}^{(-1)}\right)q_{0}^{(0)(j)}
                   + \sum\limits_{j=n-d+1}^{n}\left(\sum\limits_{\ell=n-d}^{2n-d-j}\beta_{n-\ell}^{(0)}a_{j+d-n+\ell}^{(-1)}\right)q_{0}^{(0)(j)};
                   \vspace{.3cm} \\
      u_{i}^{(-1)} &= \sum\limits_{j=i-1}^{n-d}\left(\sum\limits_{\ell=0}^{d}\beta_{\ell}^{(0)}a_{d+j-\ell}^{(-1)}\right)q_{i-1}^{(0)(j)} + \sum\limits_{j=n-d+1}^{n}\left(\sum\limits_{\ell=n-d}^{2n-d-j}\beta_{n-\ell}^{(0)}a_{j+d-n+\ell}^{(-1)}\right)q_{i-1}^{(0)(j)};\\
                   &\hspace{.1cm} \text{for} \hspace{.3cm} i=2,\ldots,n-d+1; \vspace{.3cm} \\
      u_{i}^{(-1)} &= \sum\limits_{j=i-1}^{n}\left(\sum\limits_{\ell=n-d}^{2n-d-j}\beta_{n-\ell}^{(0)}a_{j+d-n+\ell}^{(-1)}\right)q_{i-1}^{(0)(j)};
                   \hspace{.3cm} \text{for} \hspace{.3cm} i=n-d+2,\ldots,n+1. \\
    \end{split}
\end{equation*}
When $d=n$, the first and third equality shown above should be considered, and the second
sum of $u_{1}^{(-1)}$ is assumed be zero. The structure of the vector $\vec{u}^{(\mathcal{K}-1)}$ is similar to that of $\vec{u}^{(-1)}$ ($\gamma_{\centerdot}^{(\mathcal{K}-2)}$ should be used instead of $\beta_{\centerdot}^{(0)}$, coefficients $a_{\centerdot}^{(-1)}$
are substituted by $a_{\centerdot}^{(\mathcal{K}-1)}$, $q_{\centerdot}^{(0)(\centerdot)}$ is replaced by $q_{\centerdot}^{(\mathcal{K}-2)(\centerdot)}$, and the rest remain the same). For $i=1,\ldots,d$, the vector components of $\vec{v}^{(-1)}$ and $\vec{v}^{(\mathcal{K}-1)}$ are given by
\begin{equation*}
    -\sum\limits_{j=0}^{i-1}\beta_{j}^{(0)}a_{i-1-j}^{(-1)} \hspace{.3cm} \text{and} \hspace{.3cm} -\sum\limits_{j=0}^{i-1}\gamma_{j}^{(\mathcal{K}-2)}a_{i-1-j}^{(\mathcal{K}-1)},
\end{equation*}
respectively. Finally, we have $\vec{w}=\left(\sum\limits_{i=0}^{n}a_{i}^{(-1)},0,\ldots,0,a_{0}^{(\mathcal{K}-1)}\right)^{T}.$

\begin{flushleft}
\textsc{Case} 2: $d=d_{a}=d_{b}=n$ and $d_{c}<n$.
\end{flushleft}
So in this case, we have that the degree of the polynomial $c_{k}(s)$ in (\ref{ProblemaPerturbado}) is less than $n$.
From the development performed in Case 1, a square system of linear equations is obtained by considering
$c_{k}(s)=\sum\limits_{i=0}^{n}\gamma_{i}^{(k)}s^{i}$ with $\gamma_{d_{c}+1}^{(k)}=\cdots=\gamma_{n}^{(k)}=0$.

\begin{flushleft}
\textsc{Case} 3: $d=d_{a}=d_{c}=n$ and $d_{b}<n$.
\end{flushleft}
The argument is analogous to Case 2, but now applied to the polynomial $b_{k}(s)$ in (\ref{ProblemaPerturbado}) (i.e., $\beta_{d_{b}+1}^{(k)}=\cdots=\beta_{n}^{(k)}=0$).

\begin{flushleft}
\textsc{Case} 4: $d=d_{a}$, $d_{b}<n$ and $d_{c}<n$.
\end{flushleft}
We define the coefficients $\beta_{d_{b}+1}^{(k)}=\cdots=\beta_{n}^{(k)}=0$ and $\gamma_{d_{c}+1}^{(k)}=\cdots=\gamma_{n}^{(k)}=0$,
 and apply Case 1. $\Box $ \\[2mm]

\noindent{\bf Remarks.}
\begin{itemize}
\item In summary, the solution of problem~(\ref{ProblemaPerturbado}) under the hypothesis of Theorem~\ref{Teo-SolEDperturbNoAut}
is, for each step $k$, the polynomial (\ref{XsolEDperturb}), whose coefficients are obtained from solving the system of
linear equations built in the proof of Theorem~\ref{Teo-SolEDperturbNoAut}. In these circumstances, and by returning the change of scale,
the polynomial solution~(\ref{XsolEDperturb}),
for each $k=0,\ldots,\mathcal{K}-2$, will be an approximation to the exact solution of the problem~(\ref{MTFDE})-(\ref{CB2}),
since the perturbed non-autonomous differential equation~(\ref{ProblemaPerturbado}) is an approximate representation
of equation (\ref{MTFDE}).
\item We consider that the degrees of the polynomials $a_k$, $b_k$, and $c_k$ are less than or equal to $n$ and hence $d\leq n$,
because these conditions are sufficient for the purposes of numerical experimentation. The fact that we assume that $d=d_a$
 (instead of $d=d_b$ or $ =d_c$) corresponds to the fact that the number $d_a$ is important in the definition of the differential
 operator and in the construction of the canonical polynomials.
\item If $d=0$, then $d_{a}=d_{b}=d_{c}=0$ and the differential equation in (\ref{ProbDiscret}) will have constant coefficients
$a$, $b$, and $c$, as defined in (\ref{MTFDE}) for the autonomous case. The differential operator is the same for all steps
$k=0,\ldots,\mathcal{K}-2$, and its definition follows from (\ref{OperadorDif}) with $a_{k}(s)=a$, and changing the notation
from $D^{(k)}$ to $D$. Thus, from Proposition~\ref{Pro-PolCanInd} it follows that the differential operator, for autonomous case,
has no undefined canonical polynomials. The formula for generating the undefined polynomial in this case is
\begin{equation*}
    Q_{m}(s)=-\frac{m!}{a}\sum\limits_{i=0}^{m}\frac{1}{a^{m-i}i!}s^{i}, \hspace{.4cm} m\geq 0
\end{equation*}
(see \cite{DaSEsc:11}). The exact solution of the autonomous perturbed problem will be a polynomial of degree $n$,
whose coefficients and $\tau$-parameters of perturbation term (see (\ref{TermErrorProbNoAut}) with $d=0$)
are obtained by solving a system of linear equations, which can be constructed analogously as in the proof of Theorem~\ref{Teo-SolEDperturbNoAut}. Thus, the associated matrix is as defined in (\ref{matrizAprobNoAut}),
but without considering the rows containing the submatrices $R_{\gamma}^{(k)}$ for $k=0,\ldots,\mathcal{K}-3$,
$R_{\beta}^{(k)}$ for $k=1,\ldots,\mathcal{K}-2$, $\widehat{R}_{C}$ and $R_{C}$; the matrix order is
$[(n+1)(\mathcal{K}-1)+\mathcal{K}]\times [(n+1)(\mathcal{K}-1)+\mathcal{K}]$. The vector of independent terms
is formed in this case only by the vectors $\vec{u}^{(-1)}$, $\vec{u}^{(\mathcal{K}-1)}$, and $\vec{w}$,
with null vectors separating $\vec{u}^{(-1)}$ and $\vec{u}^{(\mathcal{K}-1)}$ as can be seen in (\ref{vectorbprobNoAut}).
\end{itemize}

\section{Numerical experiments}
The experiments shown in this section are defined from the following proposition, which provides us with a family
of non-autonomous MFDEs with analytical solution.
\begin{pro}
   Consider the MFDE (\ref{MTFDE}) and let $F$ be a primitive of the function $a(t)$.
   If $b(t)=-e^{F(t+1)}$ and $c(t)=e^{F(t-1)}$, then the analytical solution of non-autonomous MFDE~(\ref{MTFDE})
   is given by  $x(t)=e^{F(t)}$.
\end{pro}
\noindent{\it \textbf{Proof.}}
It is easy to verify that $x(t)=e^{F(t)}$  satisfies the non-autonomous MFDE~(\ref{MTFDE}).  $\Box$. \\[2mm]

Here, we apply the segmented Tau method obtained in the previous section to numerically solve each of the
non-autonomous problems presented hereafter. In the following experiments the boundary conditions, $\psi_{1}(t)$ and $\psi_{2}(t)$,
are defined as the analytical solution of the problem restricted to the intervals $[-1,0]$ and
$(\mathcal{K}-1,\mathcal{K}]$, respectively. We estimate the error in the infinity norm between the
numerical and analytical solutions of each problem on the interval $(0,\mathcal{K}-1]$. We consider,
on each subinterval $(k,k+1]$, $k=0,\ldots,\mathcal{K}-2$, 128 equally spaced nodes. 
All the test cases were executed using the MATLAB 7.8.0 language on an Intel Core i5-460M processor. \\[2mm]

\noindent{\bf Experiment~1.} Let us consider the following  MFDE
\begin{equation*}
    \dot{x}(t)=mx(t)-e^{m(t+1)}x(t-1)+e^{m(t-1)}x(t+1),
\end{equation*}
which has the analytical solution $x(t)=e^{mt}$, $m\in\mathbb{R}-\{0\}$.
This example is closely somehow related to numerical experiments proposed in
\cite{LiTeFoLu:10}, \cite{LTFL:10}.  \\

This non-autonomous MFDE presents a polynomial coefficient of degree $d_{a}=0$, which will form part of the
definition of the associated differential operator. This polynomial is important in the construction of the canonical
polynomials. Hence, we consider $d=d_{a}=0$. This implies that for the coefficients $b(t)=-e^{m(t+1)}$ and $c(t)=e^{m(t-1)}$
we work with a polynomial of degree zero as their estimate. Consequently, we approximate the non-autonomous MFDE by an
autonomous MFDE and estimate the solution of the latter as described in the remark at the end of Section~\ref{sectsolpert}.
\begin{table}[h]\label{Tab-Ej1}
\caption{The $\infty$-norm errors on $(0,\mathcal{K}-1]$ for $m=0.7$ and $2$.} \label{tabla1}
\begin{center}
\footnotesize{\begin{tabular}{cccccccc}
  \hline & \multicolumn{3}{l}{$m=0.7$} & & \multicolumn{3}{l}{$m=2$} \\ \cline{2-4} \cline{6-8}
  \vspace{-.3cm} & & & & & & & \\
  $n$ & $\mathcal{K}=3$ & $\mathcal{K}=5$ & $\mathcal{K}=10$ & & $\mathcal{K}=3$ & $\mathcal{K}=4$ & $\mathcal{K}=5$ \\
  \hline
  \vspace{-.3cm} & & & & & & & \\
  7 & 9.129$\times 10^{-10}$ & 4.153$\times 10^{-09}$ & 2.957$\times 10^{-07}$ & & 2.763$\times 10^{-04}$ & 1.954$\times 10^{-04}$ & 1.387$\times 10^{-02}$ \\
  8 & 1.287$\times 10^{-10}$ & 8.270$\times 10^{-11}$ & 1.123$\times 10^{-07}$ & & 7.001$\times 10^{-06}$ & 2.405$\times 10^{-04}$ & 3.430$\times 10^{-04}$ \\
  9 & 3.286$\times 10^{-13}$ & 1.585$\times 10^{-12}$ & 3.269$\times 10^{-08}$ & & 2.295$\times 10^{-07}$ & 5.517$\times 10^{-07}$ & 1.141$\times 10^{-05}$ \\
  10 & 3.331$\times 10^{-13}$ & 3.519$\times 10^{-12}$ & 7.774$\times 10^{-09}$ & & 3.672$\times 10^{-09}$ & 1.233$\times 10^{-06}$ & 5.800$\times 10^{-07}$ \\
  11 & 1.414$\times 10^{-13}$ & 2.988$\times 10^{-12}$ & 1.037$\times 10^{-08}$ & & 2.196$\times 10^{-10}$ & 7.155$\times 10^{-09}$ & 7.205$\times 10^{-07}$ \\
  12 & 3.100$\times 10^{-13}$ & 2.732$\times 10^{-11}$ & 4.548$\times 10^{-08}$ & & 9.301$\times 10^{-12}$ & 4.786$\times 10^{-09}$ & 1.075$\times 10^{-06}$ \\
  \hline
\end{tabular}}\normalsize
\end{center}
\end{table}

In Table~\ref{tabla1} the infinity-norm of the error on $(0,\mathcal{K}-1]$ is computed for $m=0.7$ with
$\mathcal{K}=3, 5, 10$ and $m=2$ with $\mathcal{K}=3, 4, 5$. Piecewise polynomial solutions
of degrees $n=7, 8, 9, 10, 11, 12$ were constructed, and satisfactory results were obtained.
Note also that when $\mathcal{K}$ becomes larger, the errors increase slightly. In the case that
$m=2$ and $\mathcal{K}=5$, the errors are greater than those shown in $m=0.7$, for the same value of $\mathcal{K}$.
The function $x(t)=e^{mt}$, with $m=2$, grows much faster than when $m=0.7$, which makes it more difficult to find
polynomial approximations to the exponential function with $m=2$ than with $m=0.7$, on unit subintervals. \\[2mm]

\noindent{\bf Experiment~2.}
Let us consider on $[0,\mathcal{K}-1]$ the following coefficients of equation~(\ref{MTFDE})
\begin{equation*}
    a(t)=\frac{3t^{2}-2t+1}{t^{3}-t^{2}+t+5}, \hspace{.5cm} b(t)=-t^{3}-2t^{2}-2t-6, \hspace{.5cm} \text{and} \hspace{.5cm} c(t)=t^{3}-4t^{2}+6t+2.
\end{equation*}
The corresponding analytical solution is $x(t)=t^{3}-t^{2}+t+5$.
\begin{table}[h]
    \caption{Segmented Tau approximation, $\infty$-norm error and perturbed terms on $(k,k+1]$ for $k=0,1$.} \label{tabla2}
    \begin{center}
        \footnotesize{\begin{tabular}{cccc}
          \hline
          Interval & Tau approximation & $\infty$-norm error & Perturbed term \\
          \hline
          \vspace{-.3cm} & & \\
          $(0,1]$ & $1.00t^{3}-0.97t^{2}+0.99t+5.00$ & 2.459$\times 10^{-02}$ & $\left(0.02t^{4}-0.12t^{3}+0.17t^{2}-0.13t-0.04\right)T_{2}^{*}(t)$ \\
          $(1,2]$ & $0.97t^{3}-0.71t^{2}+0.30t+5.45$ & 2.882$\times 10^{-02}$ & $\left(0.005t^{3}+0.002t^{2}-0.005t-0.002\right)T_{3}^{*}(t-1)$ \\
          \hline
        \end{tabular}}\normalsize
    \end{center}
\end{table}

Our approach was applied to $\mathcal{K}=3$ with $d=3$ and $n=3$, in order to estimate the polynomial solution in steps
$k=0$ and $k=1$. In Table~\ref{tabla2}, the numerical solution obtained (i.e., the ``Tau approximation") on each of the
subintervals, $(0,1]$ and $(1,2]$, is shown; in addition, we can compare the coefficients of the piecewise polynomial
solution with those of the analytical solution. Furthermore, the $\infty$-norm error and perturbed terms in each
step are shown in the third and fourth columns of Table 2, respectively; also,
$\tau$-parameters obtained by applying the method are shown in the last column of Table~\ref{tabla2}. \\[2mm]

\noindent{\bf Experiment~3.}
Let us define the coefficients of the MFDE~(\ref{MTFDE}) by
\begin{equation*}
    a(t)=\frac{\cos(t)-e^{-t}}{\sin(t)+e^{-t}+2}, \hspace{.5cm} b(t)=-\sin(t+1)-e^{-t-1}-2, \hspace{.5cm} \text{and} \hspace{.5cm} c(t)=\sin(t-1)+e^{-t+1}+2.
\end{equation*}
The analytical solution is $x(t)=\sin(t)+e^{-t}+2$.
\begin{table}[h]
    \caption{The $\infty$-norm errors on $(k,k+1]$ for $k=0,1,2,3$ with $\mathcal{K}=5$ and $d=8$.} \label{tabla3}
    \begin{center}
    \footnotesize{\begin{tabular}{ccccc}
      \hline
      $n$ & $(0,1]$ & (1,2] & (2,3] & (3,4] \\
      \hline
      \vspace{-.3cm} & & & & \\
      8 & 3.026$\times 10^{-05}$ & 8.949$\times 10^{-05}$ & 2.721$\times 10^{-04}$ & 4.876$\times 10^{-06}$ \\
      9 & 2.204$\times 10^{-04}$ & 7.235$\times 10^{-05}$ & 2.453$\times 10^{-04}$ & 1.154$\times 10^{-05}$ \\
      10 & 1.976$\times 10^{-04}$ & 1.872$\times 10^{-04}$ & 1.058$\times 10^{-04}$ & 1.039$\times 10^{-04}$ \\
      11 & 2.499$\times 10^{-06}$ & 2.236$\times 10^{-05}$ & 2.553$\times 10^{-05}$ & 1.217$\times 10^{-05}$ \\
      12 & 1.243$\times 10^{-05}$ & 1.252$\times 10^{-04}$ & 3.926$\times 10^{-05}$ & 9.230$\times 10^{-05}$ \\
      \hline
    \end{tabular}}\normalsize
    \end{center}
\end{table}

$\infty$-norm errors on $(k,k+1]$, $k=0,...,\mathcal{K}-2$ for $\mathcal{K}=5$, are shown in Table~\ref{tabla3}.
Polynomial approximations of degrees $n=8,9,10,11,12$ were computed with $d=8$. We note the
$\infty$-norm error on the entire interval $(0,4]$ matches the $\infty$-norm error obtained in $(2,3]$ for $n=8,9,11$.
While, if $n=10$, the error on $(0,4]$ is equal to that obtained on $(0,1]$ and when $n=12$, the $\infty$-norm error
on $(0,4]$ is reached in the $(1,2]$ subinterval. In this experiment, the
approximation obtained by applying the segmented Tau method provided us satisfactory numerical results. \\[2mm]

\noindent{\bf Experiment~4.}
Let $a(t)=-\frac{1}{2}(t+2)^{-3/2}$, $b(t)=-e^{\frac{1}{\sqrt{t+3}}}$, and $c(t)=e^{\frac{1}{\sqrt{t+1}}}$ for $t\in [0,\mathcal{K}-1]$.
 This problem has the analytical solution $x(t)=e^{\frac{1}{\sqrt{t+2}}}$.
 \begin{figure}[h]
  \includegraphics[height=3.0in,width=6.6in]{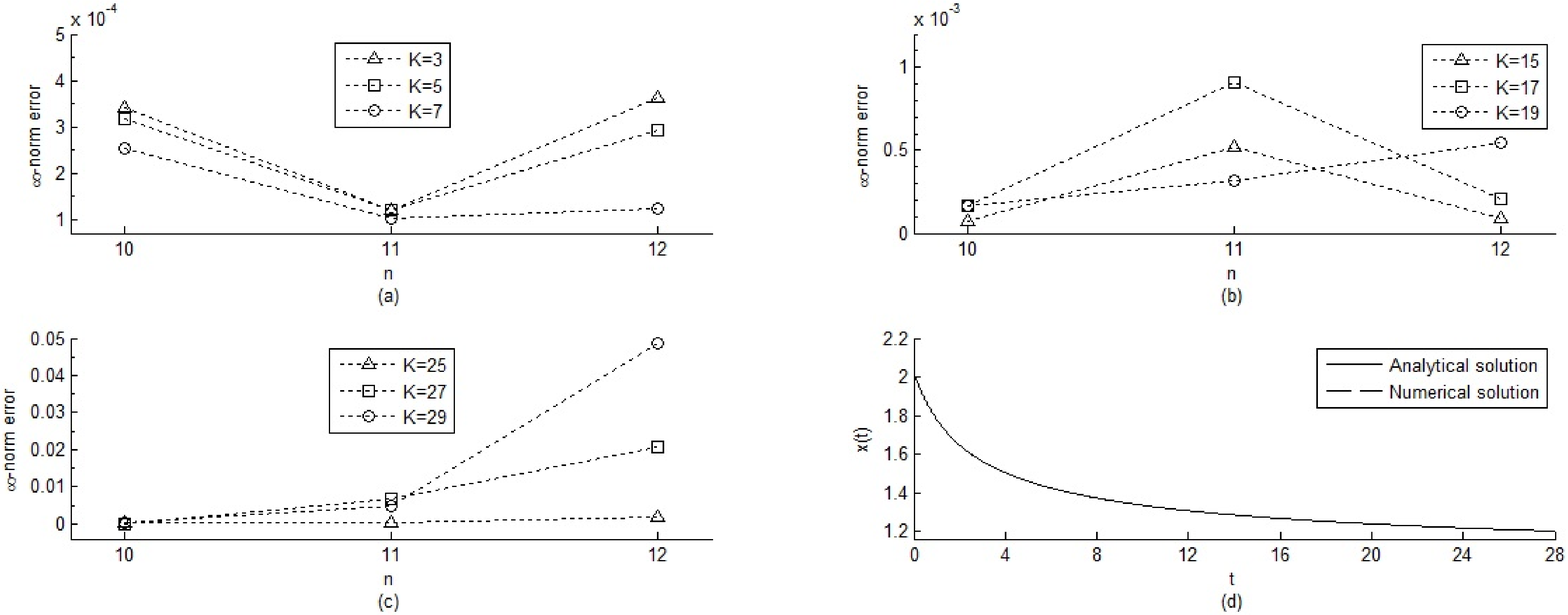}\\
  \caption{The $\infty$-norm errors for $n=10,11,12$ and graphs of curves on $(0,\mathcal{K}-1]$, (a) $\mathcal{K}=3,5,7$, (b) $\mathcal{K}=15,17,19$, (c) $\mathcal{K}=25,27,29$, (d) numerical and analytical solutions with $\mathcal{K}=29$.}\label{Fig-Ej4}
\end{figure}

We apply our approach with $d=10$. The infinity-norm errors between the numerical and analytical solutions on $(0,\mathcal{K}-1]$
were generated for different values of $\mathcal{K}$, and polynomial degrees $n=10,11,12$. In Figures~\ref{Fig-Ej4}-(a)-(b)-(c)
these errors are presented graphically. Three disjoint sets for values of $\mathcal{K}$: $\{3,5,7\}$, $\{15,17,19\}$, and $\{25,27,29\}$
were considered, which represent small, intermediate and large values of $\mathcal{K}$, respectively.
For a fixed $n$, we observe that if $\mathcal{K}\in\{25,27,29\}$, the errors are larger than when $\mathcal{K}\in\{15,17,19\}$ which,
in turn, are higher than for $\mathcal{K}\in\{3,5,7\}$.
However, for a fixed $n$, if we compare the errors on $(0,\mathcal{K}-1]$ for the three values of $\mathcal{K}$, in any of
the disjoint sets considered, we found, in general, no relationship between the error growth and the increasing values of $\mathcal{K}$
(this property is clear, for example, from Figures~\ref{Fig-Ej4}-(b)-(c) for $n = 12$).
The fact that the errors increase slightly by varying $n$, for a fixed $\mathcal{K}$, can be observed in Figure~\ref{Fig-Ej4}-(b)
with $\mathcal{K}=19$, and Figure~\ref{Fig-Ej4}-(c) for $\mathcal{K}=25,27,29$.
Figure~\ref{Fig-Ej4}-(d) shows the numerical and analytical solutions on $(0,28]$ (both graphs overlapped)
producing an infinity-norm error of 3.225$\times 10^{-04}$. In this case, the numerical solution was obtained for $n = 10$. \\[2mm]

\noindent{\bf Experiment~5.}
Let $F(t)=\ln(V(t))$, where
\begin{equation*}
V(t)=V_{p}\sin(2\pi f_{p}t)+m\frac{V_{p}}{2}\cos(2\pi(f_{p}-f_{m})t)-m\frac{V_{p}}{2}\cos(2\pi (f_{p}+f_{m})t)+\pi,
\end{equation*}
with $V_{p}=1$, $f_{p}=\frac{3}{10\pi}$, $m=\frac{1}{2}$, and $f_{m}=\frac{1}{20\pi}$.
The coefficients of the non-autonomous MFDE (1) for this example are given by
\begin{equation*}
    a(t)=F'(t), \hspace{.5cm} b(t)=-e^{F(t+1)}, \hspace{.5cm} \text{and} \hspace{.5cm} c(t)=e^{F(t-1)},
\end{equation*}
and the exact solution is $x(t)=V(t)$.\\
\begin{figure}[h]
  \includegraphics[height=3.2in,width=6.7in]{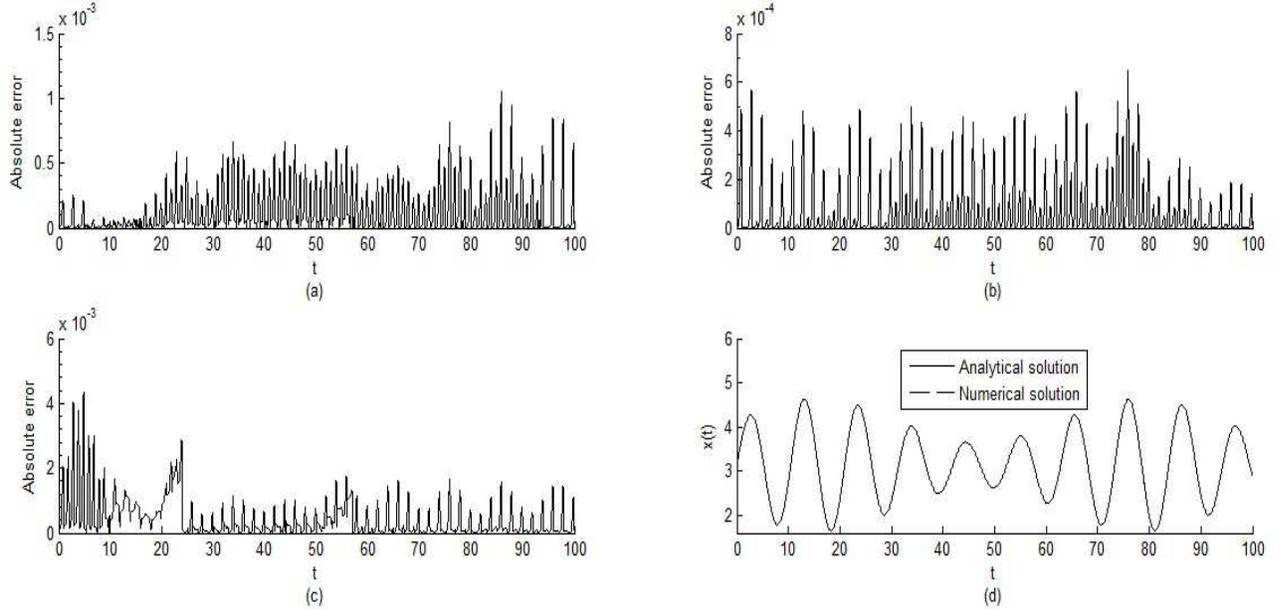}\\
  \caption{Absolute error and graphs of curves on $(0,\mathcal{K}-1]$ with $\mathcal{K}=101$ and $d=6$, (a) $n=6$, (b) $n=7$, (c) $n=8$, (d) numerical and analytical solutions.}\label{Fig-Ej5}
\end{figure}

The analytical solution has an oscillating graph with non-constant amplitude. We have considered the estimation of the approximate
solution over a domain of length 100, corresponding to $\mathcal{K}=101$. MFDE coefficients were approximated by polynomials
 of degree 6 ($=d=d_{a}=d_{b}=d_{c}$). In Figure~\ref{Fig-Ej5}-(a)-(b)-(c), we see the behavior of the absolute error of the
 numerical solution on $(0,100]$, where $n=6,7,8$. Figure~\ref{Fig-Ej5}-(d) shows the graphs of the analytical and numerical
 solutions on $(0,100]$; the numerical solution was generated with $n=7$ producing an infinity-norm error equal to
 $6.502\times 10^{-04}$.

\section{Final remarks}
In this paper a new approach to estimate the solution of a non-autonomous linear functional differential
equation of mixed type posed as a boundary value problem was presented. We adapt the segmented Tau method to the
characteristics of the problem and generate numerical solutions defined by a piecewise polynomial function.
After discretization of the problem under consideration, we define the associated linear differential operator and
present several results regarding the set of canonical polynomials associated with this type of differential
operator.
We demonstrate the existence of a system of linear equations whose solution provides the coefficients of a
polynomial solution of perturbed non-autonomous problem with polynomial coefficients and $\tau$-parameters defined
in the perturbation term.

In addition, a family of non-autonomous linear functional differential equations of mixed type with analytical
solution was provided; this result was used to formulate examples in order to show the versatility of our numerical approach.

It is worth noting, that in one of the examples treated (Experiment 1) the estimation of the solution of a non-autonomous
problem by the numerical solution of an autonomous linear functional differential equation was illustrated,
obtaining favorable results. We also show examples where numerical results on large length intervals were analyzed (e.g.,
in Experiment~5 the interval $(0,100]$ was considered).

The satisfactory results obtained in the experiments motivates us to extend in the near future this approach
 to the problem of nonlinear functional differential equations of mixed type. \\[3mm]

\noindent {\bf Acknowledgements:}
The first author was partially supported by the Consejo de Desarrollo Cient\'{\i}fico y Human\'{\i}stico (CDCH) at UCV.
The second author was partially supported by the Decanato de Investigaci\'on y Desarrollo (DID) at USB.

\bibliography{Driver}
\end{document}